\documentclass[final]{article}
\usepackage{epsfig}
\usepackage{amsmath}
\usepackage{amssymb}
\usepackage{amscd}
\usepackage{amsthm}

\begin{document}
\newtheorem*{them}{Theorem}
\newtheorem{thm}{Theorem}[section]
\newtheorem{claim}[thm]{Claim}
\newtheorem{cor}[thm]{Corollary}
\newtheorem{lemma}[thm]{Lemma}
\newtheorem{prop}[thm]{Proposition}

\theoremstyle{definition}
\newtheorem{defn}[thm]{Definition}
\newtheorem{ex}[thm]{Exercise}
\newtheorem{prob}[thm]{Problem}

\title{The Corona Theorem on the Complement of Certain Square Cantor Sets}
\author{Jon Handy   \\ University of California, Los Angeles}
\date{}
\maketitle

\begin{abstract}
Let $K$ be a square Cantor set, i.e. the Cartesian product $K=E\times E$ of two linear Cantor sets.  Let $\delta_n$ denote the proportion of the intervals removed in the $n$th stage of the construction of $E$.  It is shown that if $\delta_n=o(\frac1{\log\log n})$ then the corona theorem holds on the domain $\Omega=\mathbb C^\ast\setminus K$.
\end{abstract}

\newcommand{\diam}{{\rm{diam}}}
\newcommand{\Area}{\rm{Area}}
\newcommand{\dist}{\rm{dist}}
\newcommand{\re}{\rm{Re}\,}
\newcommand{\osc}{\rm{osc}}
\newcommand{\supp}{\rm{supp}\,}
\newcommand{\Cp}{\rm{Cap}\,}

\section{Introduction}

Denote by $H^\infty(\Omega)$ the collection of bounded analytic functions on a plane domain $\Omega$.  Let us denote by $\mathcal M=\mathcal M(H^\infty(\Omega))$ the set of multiplicative linear functionals on $H^\infty(\Omega)$, i.e. $\mathcal M$ is the maximal ideal space of $H^\infty(\Omega)$.  When $H^\infty(\Omega)$ separates the points of $\Omega$, there is a natural identification of $\Omega$ with a subset of $\Omega$ through the functionals defined by pointwise evaluations.  The corona problem for $H^\infty(\Omega)$ is to determine whether $\mathcal M$ is the closure of $\Omega$ in the Gelfand topology.

The problem can be stated more concretely as follows.  Given $\{f_j\}_1^n\in H^\infty(\Omega)$ and $\delta>0$ such that $1\ge\max_j|f_j(z)|>\delta$ for all $z\in\Omega$, do there exist $\{g_j\}_1^n\in H^\infty(\Omega)$ such that
    $$\sum_{j=1}^nf_jg_j=1?$$
We call the functions $\{f_j\}$ corona data, and the functions $\{g_j\}$ corona solutions.  It has been conjectured that this can be answered in the affirmative for any planar domain.  In the case of Riemann surfaces, a counterexample was first found by Brian Cole (see \cite{cole}), other examples being found later by D. E. Barrett and J. Diller \cite{barrett-diller}, however, this is due to a structure that in some sense makes the surface seem higher dimensional, so one might hope that the restriction to the Riemann sphere might prevent this obstacle.

This problem was first posed in the case of the unit disc $\mathbb D$ by S. Kakutani in 1941, which case was solved by L. Carleson \cite{car-disc} in 1962.  It is from this origin that the problem gets its name, as there would have been a set of maximal ideals suggestive of the sun's corona has the theorem failed in this case.  The proof was subsequently simplified by L. H\"ormander \cite{hormander} by use of a $\bar\partial$-problem, and then later in an acclaimed proof by T. Wolff (see \cite{gamelin2} or \cite{garnett-baf}).  A quite different proof of this result using techniques from several complex variables was later developed by B. Berndtsson and T. J. Ransford \cite{Berndtsson-Ransford} and Z. Slodkowski \cite{slod1}.  Besides the intrinsic interest of this result in classical function theory, the proof of this result introduced a number of tools and ideas which have proven to be of great importance in analysis.

Following the appearance of this result, it was swiftly generalized to the case of finitely connected domains, and by now a number of proofs exist for this case, e.g. \cite{alling1}, \cite{alling2}, \cite{earle-marden}, \cite{forelli}, \cite{slod2}, \cite{stout1}, \cite{stout2}, \cite{stout3}.  Although any planar domain can be exhausted by a sequence of finitely connected domains, to solve the problem on the larger domain one must control the norms of the corona solutions, $\Vert g_j\Vert_\infty$, in the approximating domains in order to take normal limits, control which is unfortunately not provided by any of the proofs just cited.

The corona problem was first solved for a class of infinitely connected domains by M. Behrens \cite{behrens1}, \cite{behrens2}.  These are ``roadrunner'' domains $\mathbb D\setminus\bigcup B_j$, where $B_j=B(c_j,r_j)$ is a disc centered at $c_j$ and radius $r_j$ such that
    $$\sum_{j=1}^\infty\frac{r_j}{|c_j|}<\infty\qquad\text{and}\qquad\left|\frac{c_{j+1}}{c_j}\right|<\lambda<1\text{ for all j}.$$
This summability restriction was improved somewhat in \cite{deeb}, \cite{deeb-wilken}.  Behrens \cite{behrens2} also proved that if the corona theorem fails for a plane domain then it fails for a domain of the form $\mathbb D\setminus\bigcup B_j$, where $\{B_j\}$ is a sequence of discs clustering only at the origin.  In this direction there is also a result of Gamelin \cite{gamelin} that the corona problem is local in that it depends only on the behavior of the domain locally about each boundary point.

The next significant progress for infinitely connected domains was again achieved by Carleson \cite{car-fatsets}, who solved the corona problem for domains having boundary $E\subset\mathbb R$ satisfying, for some $\epsilon>0$,
    $$\Lambda_1(B(x,r)\cap E)\ge\epsilon r$$
for every $x\in E$ and $r>0$, where $\Lambda_1$ denotes one-dimensional Hausdorff measure.  This proof followed an idea introduced in \cite{forelli}, constructing an explicit projection from $H^\infty(\mathbb D)$ onto $H^\infty(\Omega)$.  Following in the same vein, P. Jones and D. Marshall \cite{jones-marshall} used such projections to show that if the corona problem can be solved at the critical points of Green's function for the domain, then it can be solved for the domain, and they provided a number of sufficient conditions for this criterion.

\begin{figure}[tb]
\centering\epsfig{figure=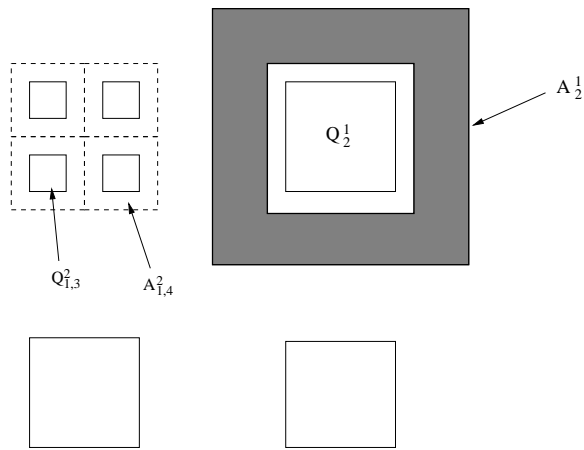, width=3in}
\caption{The various geometric constructs associated to the Cantor set:  $Q^n_J$, $V^n_J$, and $A^n_J$.}
\label{fig1}
\end{figure}

Later J. Garnett and P. Jones \cite{garnettjones} extended the result of Carleson \cite{car-fatsets}, showing that the corona theorem holds for any domain having boundary contained in $\mathbb R$.  This was later extended by C. Moore \cite{moore} to the case of domains with boundary contained in the graph of a $C^{1+\epsilon}$ function.

Due to the results in \cite{jones-marshall}, the corona problem for a domain can be solved if the critical points of Green's function (for a fixed base point) form an interpolating sequence for $H^\infty(\Omega)$.  However, it was shown by M. Gonzalez \cite{gonzalez} that for a large class of domains, conditions necessary for the critical points to be an interpolating sequence might fail.  Peter Jones (unpublished) has given an explicit example of a class of domains where the critical point of greens function fail to form an interpolating sequence, namely the complements of certain square Cantor sets (which we describe more explicitly below).

Throughout the remainder of this work, $\Omega$ will denote the complement of a square Cantor set $K$.  To fix notation, we define $K$ explicitly as $K=\bigcap_{n}K^n$, where the $K^n$ are defined inductively as follows.  Fix $\{\lambda_n\}_\mathbb N\in(\frac14,\frac12)$, which we assume satisfies $\lambda_n\le\lambda_{n+1}$ for simplicity.  (These conditions also ensure that $H^\infty(\Omega)$ is nontrivial.)  Put $K^0=[0,1]^2$.  At stage $n$, we set $K^n=\bigcup_{|J|=n}Q^n_J$ where $Q^n_J$ are squares with sides parallel to the axes and side length $\sigma_n=\prod^n_{k=1}\lambda_k$, and $J$ is a multi-index of length $|J|=n$ on letters $\{1,2,3,4\}$.  At stage $n+1$, we construct squares $Q^{n+1}_{J,j}\subset Q^n_J$, $j\in\{1,2,3,4\}$, of side length $\sigma_{n+1}$ with sides parallel to the axes such that each $Q^{n+1}_{J,j}$, $1\le j\le4$, contains a corner of $Q^n_J$.  We define $K^{n+1}=\bigcup_{|J|=n}\bigcup_{j=1}^4Q^{n+1}_{J,j}$.

Some auxiliary definitions are also useful in this setting.  The first of these is the quantity $\delta_n=1-2\lambda_n$, which represents the normalized gap width between squares of the $n$th generation having common parent.  It is useful also to introduce the thickened squares $V^n_J:=(1+\frac{\delta_n}{2\lambda_n})Q^n_J$, and from these the ``square annuli'' $A^n_J$ given by $A^n_J=\overline{V^n_J\setminus\bigcup_{j=1}^4V^{n+1}_{J,j}}$.  We note that our constants were chosen so that $\bigcup_1^4V^{n+1}_{J,j}$ is a single square concentric with and containing $Q^n_J$.  That $A^n_J\neq\emptyset$ follows from the assumption that $\lambda_n\le\lambda_{n+1}$.  Let us denote by $z^n_J$ the center of the square $Q^n_J$.  (See Figure \ref{fig1}.)

For the harmonic measures below, we will generally make the abbreviations $\omega(\cdot):=\omega(\infty,\cdot,\Omega)$ and $\omega(\cdot,z):=\omega(z,\cdot,\Omega)$ or $\omega_z(\cdot):=\omega(z,\cdot,\Omega)$ when there is no risk of confusion.

We briefly sketch the proof of the result of Jones.  Let us denote by $\{z_j\}$ the critical points of Green's function $g(z)=g(z,\infty)$ for the domain $\Omega$.  For $\{z_j\}$ to be an interpolating sequence, it is necessary that the sum $\sum g(z_j)$ be finite.  Roughly speaking, each square annulus $A^n_J$ contains one critical point, and $g(z)$ is of the same size as $\omega(Q^n_J)$ (up to a constant factor) for $z\in A^n_J$ when there are $0<a<b<\frac12$ such that $a\le\lambda_n\le b$ for every $n$.  Thus by Harnack's inequality the convergence of $\sum g(z_j)$ is equivalent to the convergence of $\sum_{n,J}\omega(Q^n_J)$.  But $\sum_{n,J}\omega(Q^n_J)=\sum_n \omega(K)=\infty$, so $\{z_j\}$ cannot be an interpolating sequence.

If one takes slightly more care in comparing harmonic measure to Green's function, this can be extended to show that $\{z^n_J\}_{n,J}$ is not an interpolating sequence when $\delta_n\le\frac{c_0}{\log n}$ for large $n$, where $c_0$ is some absolute constant.

In the other direction, it has been shown by Jones (unpublished) that the critical points of Green's function form an interpolating sequence in the case that $K$ has positive area.  (This occurs iff $\sum\delta_n<\infty$.)  This can be done via the techniques of \cite{jones-marshall} using harmonic measure estimates (see Theorems 4.1 and 4.5 of that paper), or, alternatively, by solving a $\bar\partial$-problem.

In the present work, we prove the following result.

\begin{thm}\label{corthmdim2} If $K=K(\{\lambda_n\})$ is a square Cantor set with $\delta_n=o\left(\frac1{\log\log n}\right)$ and $\{f_\mu\}^M_{\mu=1}\in H^\infty(\Omega)$ satisfy
    $$0<\eta\le\max_{1\le\mu\le M}|f_\mu(z)|\le1,$$
then there are $\{g_\mu\}_{\mu=1}^M\in H^\infty(\Omega)$ such that
    $$\sum_{\mu=1}^Mf_\mu(z)g_\mu(z)=1,\qquad z\in\Omega.$$\\
\end{thm}

This extends the result past the regime where the critical points of Green's function form an interpolating sequence.

The basic outline of the proof is as follows.  We first break the domain into a number of simply connected domains $T^n_J$ which, roughly speaking, are the cross-shaped domains $Q^n_J\setminus\bigcup_{j=1}^4Q^{n+1}_{J,j}$.  The fundamental idea is that found in \cite{garnettjones}, namely to apply Carleson's original theorem for simply connected domains and then constructively solve a $\bar\partial$-problem to obtain corona solutions for the whole domain.

To obtain the necessary cancellations, however, we require that the corona solutions on our simply connected subdomains have a certain amount of agreement on the overlaps of those domains.  To achieve this, we build these special solutions inductively, the solutions at stage $n$ obtained first by choosing solutions according to Carleson's corona theorem for simply connected domains and then solving a $\bar\partial$-problem to alter the solutions to match the neighboring solutions already constructed.  The method for solving the $\bar\partial$-problem is much like in \cite{garnettjones}, employing an interpolating sequence in the simply connected domains $T^n_J$ to build our solutions.  For our solution to have the desired special properties, however, we must choose our interpolating functions to have certain special properties, and it is here that we need the condition that $\delta_n=o(\frac1{\log\log n})$.  Section 3 is devoted to constructing these interpolating functions and solving the associated $\bar\partial$-problem in our simply connected subdomains.

Once these solutions have been constructed we paste these solutions together by solving another $\bar\partial$-problem.  In this case, the $\bar\partial$-problem is solved using the ideas of rational approximation theory (see \cite{vitushkin} or Chapter XII of\cite{gamelin-ua}).  Essentially, one solves the $\bar\partial$-problem on the intersection of two subdomains by the usual Cauchy integral representation, but in order to be able to sum these various pieces, we must add additional cancellation to each piece.  This is done by subtracting a bounded analytic function on $\Omega$ which has singular support on a portion of the Cantor set nearby and which matches derivatives of the integral at infinity.  Schwarz lemma bounds, in conjunction with the cancellations from our special corona solutions in the subdomains, then allow us to sum the terms and obtain the solution to the $\bar\partial$-problem.  This part of the proof implicitly makes use of the fact that for any $\zeta\in K$ there is $c_0>0$ such that
    $$\gamma(B(\zeta,r)\cap K)\ge c_0r$$
for $r\in(0,\diam(K)]$, where
    $$\gamma(E):=\sup\{|f'(\infty)|:f\in H^\infty(\mathbb C^\ast\setminus E),\Vert f\Vert\le1\}$$
denotes the analytic capacity of the set $E$ (see \cite{tmv}).  This parallels the thickness condition for the boundary employed in \cite{car-fatsets}.  Thus the functions introduced to obtain cancellations are polynomials in the extremal function for this problem, the Ahlfor's function, for a piece of the boundary.  Alternatively, one can take the function to be powers of the Cauchy integral of the uniform measure on an appropriate subsquare $K^n_J:=K\cap Q^n_J$ of the Cantor set (see \cite{garnett-ac}).

The present work is part of the author's Ph.D. dissertation.  The author is naturally greatly indebted to his advisor, John Garnett, for many helpful conversations and suggestions over the years, and would like to take this opportunity to express his gratitude.  The author would also like to thank Peter Jones for helpful conversations.

\section{Proof of the Theorem}

\begin{figure}[tb]
\centering\epsfig{figure=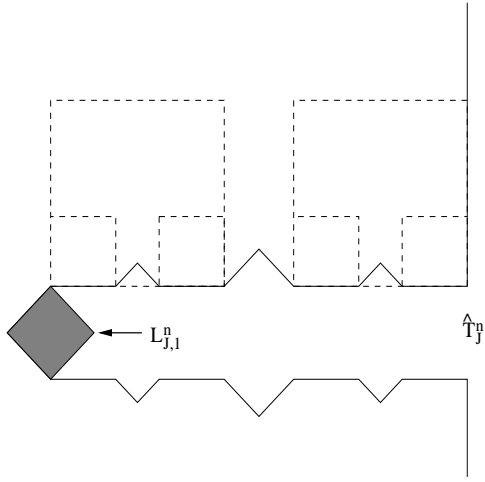,height=2.5in,width=2.5in}
\caption{One of the domains $\hat T^n_J$.}
\label{fig2}
\end{figure}

To begin the proof, let us first define simply connected, cross-shaped regions $T^n_J$ by $T^n_J:=Q^n_J\setminus\bigcup_1^4Q^{n+1}_{J,j}$.  For these regions we distinguish the four boundary segments, denoted $\{\ell^n_{J,j}\}_{j=1}^4$, which do not lie on a side of a square $Q^{n+1}_{J,j}$, $1\le j\le4$.  Let us then define lozenges
    \begin{equation}\label{lozenges}
    L^n_{J,j}=\left\{z:\left|\arg\frac{z-x_2}{x_2-x_1}\right|\vee\left|\arg\frac{z-x_1}{x_2-x_1}\right|<\alpha\right\},
    \end{equation}
where $(x_1,x_2)=\ell^n_{J,j}$ and $\alpha$ is some angle less than $\frac\pi4$ to be determined by the constructions below.  We now define regions
    $$\hat T^n_J=T^n_J\cup\bigcup\{L^m_{I,i}:\ell^m_{I,i}\subset\partial T^n_J\}.$$
(See Figure \ref{fig2}.)  We note that these regions remain simply connected, so that Carleson's corona theorem for the unit disc $\mathbb D$ provides corona solutions $\{g^{(n,J)}_\mu\}_{\mu=1}^M$ corresponding to $\{f_\mu\}_1^M$ in $\hat T^n_J$.  We note that $\sum\chi_{\hat T^n_J}\le2$.  The most important result for the present construction is embodied in the following proposition.

\begin{prop}\label{basicconstruction} Given corona data $\{f_\mu\}_1^M$ in $\Omega$ as above, if the angle $\alpha$ defining the lozenges $L^n_{J,j}$ is sufficiently small then there are corona solutions $\{g^{(n,J)}_\mu\}_1^M$ in each $\hat T^n_J$ such that $\Vert g^{(n,J)}_\mu\Vert_\infty\le C(M,\delta,K)$ and if $(n,J)$ and $(m,I)$ are indices with $\hat T^n_J\cap \hat T^m_I=L^n_{J,j}$ then for $1\le\mu\le M$,
    $$|g^{(n,J)}_\mu(z)-g_\mu^{(m,I)}(z)|\lesssim\frac1{n^3}\frac{d(z,K)}{\sigma_n\delta_n},\qquad z\in L^n_{J,j}.$$\\
\end{prop}

Assuming this proposition, let us prove the theorem.  Let $\{\phi_{(n,J)}\}$ to be a partition of unity on $\Omega$ subordinate to $\hat T^n_J$ satisfying $|\nabla\phi_{(n,J)}(z)|\lesssim d(z,K)^{-1}$ and define
    \begin{equation}\label{dbar1}
    \tilde g_\mu=\sum_{(n,J)}\phi_{(n,J)}g^{(n,J)}_\mu.
    \end{equation}
We note that the functions $\tilde g_\mu$, while not analytic, have the desired property that $\sum f_\mu\tilde g_\mu\equiv1$ on $\Omega$.  An observation due to H\"ormander \cite{Hormander}, now reduces us to solving the $\bar\partial$-problem
    $$\bar\partial a_{\mu\nu}=\tilde g_\mu\bar\partial\tilde g_\nu,\qquad a_{\mu\nu}\in L^\infty(\Omega).$$
Indeed, given such functions $a_{\mu\nu}$, if we define
    $$G_\mu=\tilde g_\mu+\sum_{\nu=1}^M(a_{\mu\nu}-a_{\nu\mu})f_\nu$$
the antisymmetry of the matrix $A=[a_{\mu\nu}]$ gives us $\sum f_\mu G_\mu\equiv1$ on $\Omega$, while (\ref{dbar1}) provides $\bar\partial G_\mu\equiv0$ on $\Omega$ for each $1\le\mu\le M$, so that $G_\mu\in H^\infty(\Omega)$.  Thus the functions $\{G_\mu\}_1^M$ provide the corona solutions sought by the theorem.

Let us therefore turn our attention to solving (\ref{dbar1}).  The immediate thought is to consider
    $$a_{\mu\nu}(z)=\frac1{2\pi i}\iint_\Omega\frac{\tilde g_\mu(\zeta)\bar\partial\tilde g_\nu(\zeta)}{\zeta-z}\,d\zeta\,d\bar\zeta,$$
but it is not clear that the integral is convergent \`a priori, so instead we begin by viewing this formally as
    $$\frac1{2\pi i}\iint_\Omega\frac{\tilde g_\mu(\zeta)\bar\partial\tilde g_\nu(\zeta)}{\zeta-z}\,d\zeta\,d\bar\zeta
        =\sum\frac1{2\pi i}\iint_{L^n_{J,j}}\frac{\tilde g_\mu(\zeta)\bar\partial\tilde g_\nu(\zeta)}{\zeta-z}\,d\zeta\,d\bar\zeta,$$
and attempt to introduce additional cancellations in each term.  We begin by estimating individual terms.  For ease of notation, we define
    $$I^n_{J,j}(z) = \frac1{2\pi i}\iint_{L^n_{J,j}}\frac{\tilde g_\mu(\zeta)\bar\partial\tilde g_\nu(\zeta)}{\zeta-z}\,d\zeta\,d\bar\zeta.$$
Let us suppose that $(m,I)$ is such that $L^n_{J,j}=\hat T^n_J\cap\hat T^m_I$.  We first note that, for $z\in\overline{L^n_{J,j}}\setminus\{x_1,x_2\}$,
    \begin{eqnarray*}
    |I^n_{J,j}(z)|
        &\lesssim& \iint_{L^n_{J,j}}\frac{|g_\nu^{(n,J)}\bar\partial\phi_{(n,J)}+g_\nu^{(m,I)}\bar\partial\phi_{(m,I)}|}{|\zeta-z|}\,dx\,dy   \\
        &\lesssim& \iint_{L^n_{J,j}}\frac{|g_\nu^{(n,J)}-g_\mu^{(m,I)}||\bar\partial\phi_{(n,J)}|}{|\zeta-z|}\,dx\,dy  \\
        &\lesssim& \frac1{\sigma_n\delta_nn^3}\iint_{L^n_{J,j}}\frac{dx\,dy}{|\zeta-z|}  \\
    |I^n_{J,j}(z)| &\lesssim& \frac1{n^3},
    \end{eqnarray*}
where we have exploited the fact that $\sum\chi_{\hat T^n_J}\le2$ in the second line and the proposition in the third.  Proceeding in a similar manner, we can achieve estimates
    $$\left|\iint_{L^n_{J,j}}\tilde g_\mu\bar\partial g_\nu\,d\zeta\,d\bar\zeta\right|\lesssim\frac{\sigma_n\delta_n}{n^3}.$$

To generate further cancellations, we employ the ``derivative matching trick'' from rational approximation theory (see \cite{vitushkin} or Chapter XII of \cite{gamelin-ua}).  To this end, we introduce functions $k^n_{J,j}$ defined as follows.  Let $\tilde K^n_{J,j}$ be a ``square'' $K^{m'}_{I'}=K\cap Q^{m'}_{I'}$ with $\tilde K^n_{J,j}\cap L^n_{J,j}\neq\emptyset$ and side length comparable to $\sigma_n\delta_n$.  The analytic capacity $\gamma(\tilde K^n_{J,j})$ is comparable to $\sigma_n\delta_n$, so if $f^n_{J,j}$ is the Ahlfors function for $\tilde K^n_{J,j}$ then, choosing (uniformly bounded) constants $c_{J,j}^n$ appropriately, and setting
    $$k^n_{J,j}(z)=\frac{c^n_{J,j}}{n^3}f^n_{J,j}(z)$$
then $\Vert k^n_{J,j}\Vert_\infty\lesssim\frac1{n^3}$ and
    $$(k^n_{J,j})'(\infty)=\frac1{2\pi i}\int_{L^n_{J,j}}\tilde g_\mu\bar\partial\tilde g_\nu\,d\zeta\,d\bar\zeta.$$
For the estimates that follow, we will employ the following form of Schwarz's lemma.

\begin{lemma}[Schwarz's Lemma] If $E$ is a compact set and $f\in H^\infty(\mathbb C^\ast\setminus E)$ has a double zero at infinity then
    $$|f(z)|\lesssim\frac{\Vert f\Vert_\infty\diam(E)^2}{d(z,E)^2}.$$\\
\end{lemma}

Let us define
    $$h^n_{J,j}(z):=I^n_J(z)-k^n_{J,j}(z).$$
Applied in the current context, the lemma yields
    \begin{equation}\label{Schwarzest}
    |h^n_{J,j}(z)|\lesssim\frac1{n^3}\left(1\wedge\frac{(\sigma_n\delta_n)^2}{d(z,\tilde K^n_{J,j}\cup L^n_{J,j})^2}\right).
    \end{equation}
Since $\bar\partial k^n_{J,j}=0$ in $\Omega$, formally we have
    $$\bar\partial\sum_{(n,J,j)}h^n_{J,j}=\tilde g_\mu\bar\partial\tilde g_\nu,$$
so it suffices to check the boundedness of the sum.

Fix $z\in\Omega$.  For each $n\in\mathbb N$ there are at most boundedly many terms for which the minimum is one in the inequality (\ref{Schwarzest}), and the summability of $\frac1{n^3}$ shows that these terms, summed over $n$, give a contribution which is controlled by the sum $\sum\frac1{n^3}$.

For remaining terms, we distinguish between the cases $\sigma_n\gtrsim d(z,K)$ and $\sigma_n\lesssim d(z,K)$.  Let $n_0\in\mathbb N$ be such that $d(z,K)\asymp \sigma_{n_0}$ (the maximum principle prevails over any $z$ with $d(z,K)\gg1$).  For $m\ge n_0$, all remaining terms are at (roughly) distance $\sigma_{n_0}$ or more away.  Inductively one can show that there are at most $4^{m-n_0+k}$ terms having singularities a distance $\sigma_{n_0-k}$ away from $z$, and so Schwarz lemma bounds the sum of these terms by a constant multiple of
    \begin{eqnarray*}
    \frac1{m^3}\sum_{k=0}^{n_0}\frac{\sigma_m^2\delta_m^2 4^{m-n_0+k}}{\sigma_{n_0-k}^2}
        &\le& \frac1{m^3}\sum_{k=0}^{n_0}\frac{4^m\sigma_m^2}{4^{n_0-k}\sigma_{n_0-k}^2}\le\frac1{m^3}\sum^{n_0}_{k=0}1  \\
        &\le& \frac{n_0}{m^3}\le\frac1{m^2}
    \end{eqnarray*}
Similarly, if $m<n_0$, there are $4^{k-m}$ terms at roughly distance $\sigma_{m-k}$, and so for these terms the sum is controlled by $m^{-3}\sum_0^m\frac{\sigma_m^24^{k-m}}{\sigma_{m-k}^2}\lesssim m^{-2}$.  Summing these bounds over $m$, we obtain $\Vert\sum_{n,J,j}h^n_{J,j}\Vert_{L^\infty(\Omega)}\lesssim1$ as desired.  This proves the theorem.\qed\\

\section{Proof of Proposition \ref{basicconstruction}}

Let us now turn to the proof of Proposition \ref{basicconstruction}.  By normal families, it suffices to perform the construction in the domains $\{\hat T^n_J:|J|=n,n\le N\}$, provided we obtain constants independent of $N$.

\begin{figure}[tb]
\centering\epsfig{figure=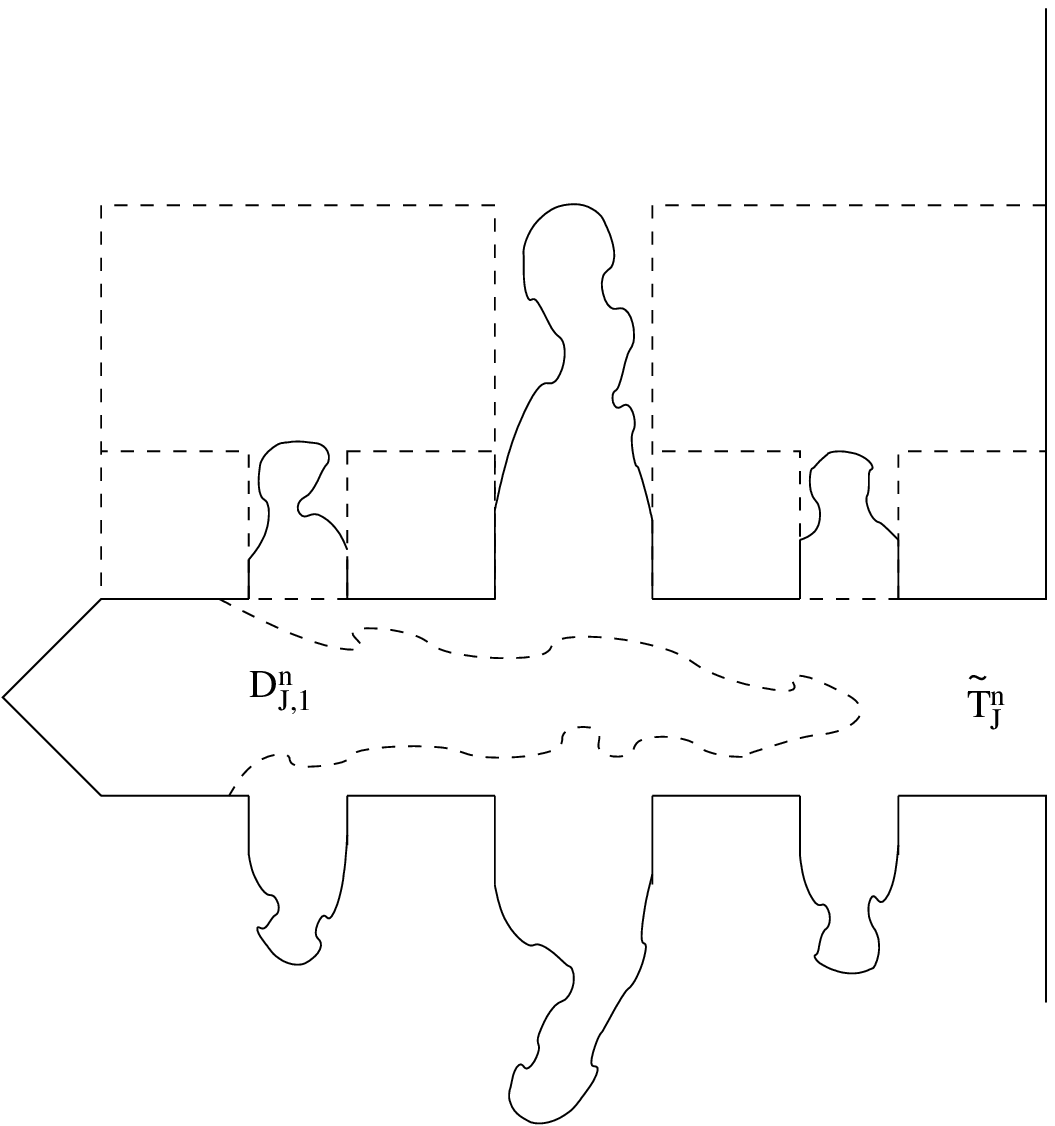,height=2.5in,width=2.5in}
\caption{One of the regions $\tilde T^n_J$.}
\label{fig3}
\end{figure}

Rather than working directly with the sets $\hat T^n_J$, we will consider slightly enlarged domains $\tilde T^n_J$, defined as follows.  We first define $\tilde T^N_J=\hat T^N_J$.  Now let $\phi^N_J$ be a conformal mapping from $\tilde T^N_J$ to $\mathbb D$ preserving the symmetries of the domain $\tilde T^N_J$.  Let $E_j=E_j(N,J)$ be the arc of $\partial\mathbb D$ defined by $E_j:=\phi^N_J(\partial L^N_{J,j}\cap\partial\hat T^N_J)$, and let $E_j^\ast=3E_j$.  Defining $\theta_N=\pi(1-\frac{c_1}{\log N})$, the constant $c_1$ to be determined later, let $\gamma_j^\ast$ be the circular arc in $\mathbb D$ with endpoints coincident with those of $E_j^\ast$ and intersecting at an angle of $\theta_N$.  Easy length-area/extremal length estimates in $T^N_J$ yield $\omega(z^N_J,\ell^N_{J,j},T^N_J)\lesssim e^{-\frac{c_0}{\delta_N}}$ for some constant $c_0$ not depending on $N$.  Employing the maximum principle and following the mapping to the disc, we find that the length of the arcs $\gamma_j^\ast$ is $o(\frac1{\log N})$ by our assumptions on the sequence $\{\delta_n\}$.  Now by elementary geometry, the disc determined by $\gamma_j^\ast$ has radius $r$ bounded by
    $$r\lesssim\frac{c_1}{\log N}e^{-\frac{c_0}{\delta_N}},$$
so that $r=o(1)$.  The arcs $\gamma_j^\ast$, $1\le j\le 4$, are therefore (uniformly) hyperbolically separated, at least for $N$ greater than some $n_0$.  Let us now define $D^N_{J,j}$ to be the simply connected domain which is the pre-image under $\phi^N_J$ of the domain bounded by $E_j^\ast\cup\gamma_j^\ast$.  We then define domains
    $$\tilde T^{N-1}_I:=\hat T^{N-1}_I\cup\bigcup\{D^N_{J,j}:L^N_{J,j}\subset\hat T^{N-1}_I\}$$
for each multi-index $I$ of length $N-1$.

Proceeding inductively, let us suppose that the domains $\tilde T^n_J$ have already been constructed for $n>m\ge n_0$, and, fixing $(m+1,J)$, let $\phi^{m+1}_J$ be a conformal mapping from $\tilde T^{m+1}_J$ preserving the symmetries of the domain.  Let $E_j=E_j(m+1,J)(\partial L^{m+1}_{J,j}\cap\partial\hat T^{m+1}_J)$, let $E_j^\ast=3E_j$, and let $\gamma_j^\ast$ be the circular arc in $\mathbb D$ with endpoints coincident with those of $E_j^\ast$ and intersecting at an angle of $\theta_{m+1}=\pi(1-\frac{c_1}{\log(m+1)})$.  As above, length-area estimates in $T^{m+1}_J$ yield $\omega(z^{m+1}_J,\ell^{m+1}_{J,j},T^{m+1}_J)\lesssim e^{-\frac{c_0}{\delta_{m+1}}}$, and thus the length of the arcs $\gamma_j^\ast$ is $o(\frac1{\log(m+1)})$ by our assumptions on the sequence $\{\delta_n\}$.  Moreover, the disc determined by $\gamma_j$ has radius $r=r_m$ bounded by
    \begin{equation}\label{radiusest}r\lesssim\frac{c_1}{\log(m+1)}e^{-\frac{c_0}{\delta_{m+1}}},\end{equation}
so that $r_m=o(1)$.  The arcs $\gamma_j^\ast$, $1\le j\le 4$, are therefore hyperbolically separated since we have taken $m\ge n_0$.  Let us now define $D^{m+1}_{J,j}$ to be the pre-image under $\phi^{m+1}_J$ of the domain bounded by $E_j^\ast\cup\gamma_j^\ast$.  We then define domains
    $$\tilde T^m_I:=\hat T^m_I\cup\bigcup\{D^n_{J,j}:n>m,L^n_{J,j}\subset\hat T^m_I\}$$
for each multi-index $I$ of length $m$.  (See Figure \ref{fig3}.)  We note that for $n\ge n_0$ the estimate (\ref{radiusest}) yields that these domains are simply connected and also that $\tilde T^n_J\cap\tilde T^n_I=\emptyset$ when $I\neq J$.

For $n<n_0$ we will take
    $$\tilde T^n_J=\hat T^n_J\cup\bigcup\{D^m_{I,j}:n\ge n_0,L^m_{I,j}\subset\hat T^n_J\}.$$

Proceeding from these definitions, we will construct our solutions from the top down, exploiting the natural generations structure of the Cantor set.  At each stage we will obtain corona solutions $\{g^{(n,J)}_\mu\}\in H^\infty(\tilde T^n_J)$, and the corona solutions specified by the proposition will simply be the restrictions of these solutions to $\hat T^n_J$.  For $n\le n_0$, we construct our corona solutions by applying the corona theorem for finitely connected domains to the domain $\Omega_{n_0}=\mathbb C^\ast\setminus[0,1]^2\cup\bigcup_{n\le n_0}\bigcup_{|J|=n}\tilde T^n_{J,j}$.  Given $n\ge n_0$, let us suppose we have already obtained the desired corona solutions $\{g^{(m,I)}_\mu\}$ in the regions $\tilde T^m_I$ for which $m<n$.  As the domains $\{\tilde T^n_J\}_{|J|=n}$ are disjoint, we may construct our corona solutions in each of those domains separately.

Fix $\tilde T^n_J$.  We note that $\tilde T^n_J\cap\left(\bigcup_{m<n}\bigcup_{|I|=m}\tilde T^m_I\right)=\bigcup_1^4D^n_{J,j}$, and so in each region $D^n_{J,j}$ there are corona solutions constructed from previous generations, which we shall denote $\{g^j_\mu\}_\mu$.  (Since $\sum_{n,J}\chi_{\tilde T^n_J}\le2$, these solutions are uniquely determined among those previously constructed.)  We now push the situation forward to the unit disc according to the map $\phi^n_J$.  Let $D_j$ denote the push-forward of the domain $D^n_{J,j}$, $F_\mu$ the push-forward of a corona datum $f_\mu$, and $G_\mu^j$ the push-forward of $g^j_\mu$.  We note that $D_j$ is a lens-shaped domain by construction.

In order to construct our corona solutions, we will, as above, reduce the problem to an appropriate $\bar\partial$-problem.  Obtaining the desired bounds in this manner will require the use of interpolating sequences, so before continuing with the main line of argument we make a detour to obtain the special interpolating functions that we require.

\begin{lemma}\label{outerfunction}  There is a function $G\in H^\infty(\mathbb D)$ of norm one with the following properties.  \begin{enumerate}
    \item There is a constant $c_2$, not depending on $n$, such that $|G(z)|\ge c_2$ for $z\in\gamma_j^\ast$.
    \item If $\tilde\alpha<\frac\pi2$ and $L_j$ is the lens domain in $\mathbb D$ bounded by $E_j$ and the circular arc $\gamma_j$ meeting the endpoints $\zeta_1(j),\zeta_2(j)$ of $E_j$ in angle $\tilde\alpha$ then
    $$|G(z)|\le\frac1{n^3}\frac{d(z,\{\zeta_1,\zeta_2\})}{\diam(L_j)}.$$
\end{enumerate}\end{lemma}
{\bf Proof: } We will obtain $G$ by multiplying a number of outer functions.  For the first factor, we define
    $$H_0(z)=\exp\left\{-\int_{\partial\mathbb D}\frac{e^{i\theta}+z}{e^{i\theta}-z}6\log n\,
        \chi_{\bigcup E_j^\ast}(\theta)\frac{d\theta}{2\pi}\right\},$$
noting that
    \begin{eqnarray*}
    |H_0(z)| &=& \exp\left\{-\int_{\partial\mathbb D}P_z(\theta)6\log n\,\chi_{\bigcup E_j^\ast}(\theta)\frac{d\theta}{2\pi}\right\}    \\
        &=& \exp\left\{-6\log n\,\omega\bigg(z,\bigcup E_j^\ast,\mathbb D\bigg)\right\},
    \end{eqnarray*}
where $P_z$ is the Poisson kernel.  Due to our choice of angle $\theta_n$, we have
    $$\omega\bigg(z,\bigcup E_j^\ast,\mathbb D\bigg)\le4\omega(z,E_k,\mathbb D)\le\frac{4c_1}{\log n},\qquad z\in\gamma_k^\ast,$$
so that $|H_0(z)|\ge e^{-24c_1}$ on each $\gamma_j^\ast$.  Also, on $L_j$ we can easily compute that $\omega(z,E_j^\ast,\mathbb D)\ge\frac12$, whereby $|H_0(z)|\le\frac1{n^3}$ on $L_j$.

For the remaining factors we first note that on the imaginary axis the functions
    $$u_\pm(z) := \int_{-\frac\pi2}^\frac\pi2\log^-|e^{\pm i\frac\pi6}-e^{i\theta}|P_{iy}(\theta)\frac{d\theta}{2\pi} $$
are bounded below by $-\log2$, and so after precomposing these functions by appropriate M\"obius transformations and exponentiating, we obtain functions $\tilde H_{j,1}$, $\tilde H_{j,2}$ of norm one such that $|\tilde H_{j,1}(z)|,|\tilde H_{j,2}(z)|\ge\frac12$ on each $\gamma_k$, $1\le k\le 4$, and such that
    \begin{eqnarray*}
    |\tilde H_{j,1}(z)| &\le& \frac{|z-\zeta_1(j)|}{\diam(L_j)},   \\
    |\tilde H_{j,2}(z)| &\le& \frac{|z-\zeta_2(j)|}{\diam(L_j)},
    \end{eqnarray*}
for $z\in L_j$.  Setting $H_j=\tilde H_{j,1}\tilde H_{j,2}$, and then $G=\prod_0^4H_j$ thus provides the desired function.\qed\\

Fixing $\beta>0$, if $\{z_k\}\in\gamma_j$ are points satisfying
    \begin{equation}\label{sepbeta}|z_k-z_\ell|\ge\beta(1-|z_k|),\qquad k\neq\ell,\end{equation}
then $\sum_k(1-|z_k|)\delta_{z_k}$, where $\delta_{z_k}$ in this case denotes the unit point mass at $z_k$, is a Carleson measure with norm depending only on $\beta$, and so $\{z_k\}$ is an interpolating sequence with constant of interpolation depending only on $\beta$ by Carleson's interpolation theorem for $H^\infty(\mathbb D)$ \cite{Carint}.  Due to the hyperbolic separation of the arcs $\gamma_j$ in $\mathbb D$, this remains true for a sequence $S=\bigcup_1^4 S_j$ with each $S_j$ a sequence in $\gamma_j$ satisfying this condition.  Then if $G(z)$ is the function provided by Lemma \ref{outerfunction}, given $\{w_j\}\in\ell^\infty$ we can find an interpolating function $f\in H^\infty(\mathbb D)$ with $f(z_j)=\frac{w_j}{G(z_j)}$.  The function $g=fG\in H^\infty(\mathbb D)$ then satisfies
    $$g(z_j)=w_j,$$
    $$\Vert g\Vert_{L^\infty(\mathbb D)}\le A'c_0,$$
    $$|g(z)|\le\frac{A'd(z,\{\zeta_1,\zeta_2\})}{n^3\diam(L_j)},\qquad z\in L_j,$$
for $L_j$ and $\zeta_1(j),\zeta_2(j)$ defined as in the lemma, where $A'$ bounds the largest constant of interpolation for a maximal sequence on $\bigcup_1^4\gamma_j$ satisfying (\ref{sepbeta}).  Fixing $\beta$ for the remainder of the proof, let $A=c_0A'$.

\begin{lemma}\label{lem2}  Let $S=\{z_n\}$ be a maximal sequence on $\bigcup_1^4\gamma_j$ satisfying
    $$|z_k-z_\ell|\ge\frac{1-|z_j|}{8A^2},\qquad k\neq\ell.$$
Then there are functions $h_j\in H^\infty(\mathbb D)$ such that
    $$h_j(z_j)=1,$$
    $$\Vert h_j\Vert_{L^\infty(\mathbb D)}\le A^2,$$
    $$\sum_j|h_j(z)|\le\frac{\log8A^2}{\log\beta^{-1}}A^2,$$
and such that
    $$\sum_j|h_j(z)|\le\frac{\log8A^2}{\log\beta^{-1}}\frac{A^2}{n^6}\left(\frac{d(z,\{\zeta_1,\zeta_2\})}{\diam(L_j)}\right)^2,\qquad z\in L_j,$$
where $L_j$ and $\zeta_1(j),\zeta_2(j)$ are defined as in the Lemma \ref{outerfunction}.
\end{lemma}
{\bf Proof: } This is a mild refinement of an argument of Varopoulos \cite{Var}.  We note that the sequence $S$ can be split into $\frac{\log8A^2}{\log\beta^{-1}}$ disjoint sequences $S_m$ such that (\ref{sepbeta}) holds.  Restricting our attention to a subsequence $S_m$, it suffices to consider the case that $S_m$ is finite, $S_m=\{z_1,\ldots,z_{n_0}\}$, as one may then employ normal families to the construction below.  Set $\omega=e^\frac{2\pi i}{n_0}$.  Employing the remarks above, by Lemma \ref{outerfunction} we may choose $f_j\in H^\infty(\mathbb D)$ such that $f_j(z_k)=\omega^{jk}$, $\Vert f_j\Vert_{L^\infty(\mathbb D)}\le A$, and $|f_j(z)|\le \frac{A}{n^3}\frac{d(z,\{\zeta_1,\zeta_2\})}{\diam(L_\ell)}$ when $z\in L_\ell$ for each $1\le\ell\le4$.  If we now define
    $$h_j(z)=\left(\frac1{n_0}\sum_{k=1}^{n_0}\omega^{-jk}f_k(z)\right)^2,$$
then
    $$h_j(z_i)=\left(\frac1{n_0}\sum^{n_0}_{k=1}\omega^{(i-j)k}\right)^2
        =\left\{\begin{array}{ll}1\qquad&\mbox{if }i=j,\\0& \mbox{else,}\end{array}\right.$$
    $$\Vert h_j\Vert_{L^\infty(\mathbb D)}\le A^2,$$
and
    $$|h_j(z)|\le \frac{A^2}{n^6}\left(\frac{d(z,\{\zeta_1,\zeta_2\})}{\diam(L_\ell)}\right)^2$$
for $z\in L_\ell$.  Also,
    \begin{eqnarray*}
    \sum^{n_0}_{j=1}|h_j(z)| &=& \frac1{n_0^2}\sum^{n_0}_{k=1}\sum_{j,\ell}\omega^{-jk}\omega^{j\ell}f_j(z)\bar f_\ell(z)    \\
        &=& \frac1{n_0^2}\sum^{n_0}_{j=1}n_0|f_j(z)|^2,
    \end{eqnarray*}
so that
    $$\sum_{j=1}^{n_0}|h_j(z)|\le A^2$$
throughout the unit disc, and
    $$\sum^{n_0}_{j=1}|h_j(z)|\le \frac{A^2}{n^6}\left(\frac{d(z,\{\zeta_1,\zeta_2\})}{\diam(L_\ell)}\right)^2,$$
for $z\in L_\ell$ as desired.  \qed\\

With these lemmas in hand we turn now to constructing our special corona solutions.  By Carleson's original corona theorem \cite{car-disc}, there are corona solutions $\{g_\mu\}$ in $\mathbb D$ with $\Vert g_\mu\Vert\le c(M,\eta)$.  To generate corona solutions close to $G_\mu^j$ on $D_j$ our first instinct is to take a partition of unity and paste these together.  In doing this we first add to the domains $D_j$ to generate overlap.  Specifically, let us define
    $$\tilde D_0 = \bigg(\mathbb D\setminus\bigcup_1^4D_j\bigg)\cup\bigcup_{z_j\in S}B\left(z_j,\frac{1-|z_j|}{4A^2}\right)$$
and
    $$\tilde D_j = D_j\cup\bigcup_{z_k\in S\cap\gamma_j}B\left(z_j,\frac{1-|z_j|}{4A^2}\right)$$
for $1\le j\le4$.  For ease of notation in what follows, we will denote the region of overlap, $\bigcup_1^4(\tilde D_j\cap\tilde D_0)$, by $U$.  Let $\psi_0,\ldots,\psi_4$ be a partition of unity subordinate to $\tilde D_0,\ldots,\tilde D_4$ satisfying $|\nabla\psi_j(z)|\lesssim (1-|z|)^{-1}$.  Our initial pasting is then
    $$\tilde g_\mu:=\psi_0g_\mu+\sum_{j=1}^4\psi_jG_\mu^j.$$

To obtain a bounded analytic solution from this we now take an approach much like that above.  In particular, we wish to find functions $a_{\mu\nu}\in L^\infty(\mathbb D)$ such that
    $$\bar\partial a_{\mu\nu}=\tilde g_\mu\bar\partial\tilde g_\nu$$
and
    $$|a_{\mu\nu}(z)|\lesssim n^{-6}\left(\frac{d(z,\{\zeta_1,\zeta_2\})}{\diam(L_j)}\right)^2,$$
for $z\in L_j$, where $L_j$ is the lens domain bounded by $E_j$ and the circular arc meeting the endpoints $\zeta_1(j),\zeta_2(j)$ of $E_j$ in angle $\tilde\alpha<\frac\pi2$.

\begin{lemma}\label{lem3} If $B\in L^\infty(U)$ and $b(z)=\frac{B(z)}{1-|z|}\chi_U$ then there is $F\in L^\infty(\mathbb D)$ such that
    $$\bar\partial F=b$$
in the sense of distributions on $\mathbb D$,
    $$\Vert F\Vert_\infty\lesssim\Vert B\Vert_\infty,$$
and
    $$|F(z)|\lesssim n^{-6}\left(\frac{d(z,\{\zeta_1,\zeta_2\})}{\diam(L_j)}\right)^2,\qquad z\in L_j,$$
where $L_j$ is the lens domain bounded by $E_j$ and the circular arc meeting the endpoints $\zeta_1(j),\zeta_2(j)$ of $E_j$ in angle $\tilde\alpha<\frac\pi2$.
\end{lemma}
{\bf Proof: } We follow an argument due to Peter Jones \cite{garnettjones}.  Let $\{h_m\}$ be the functions provided by Lemma \ref{lem2}, and let us write $U$ as the disjoint union of sets $U_m\subset B(z_m,\frac{1-|z_m|}{4A^2})$, and let us define
    $$F(\zeta) = \sum_m\frac1\pi\iint_{U_m}\frac{h_m(\zeta)}{h_m(z)}\frac{b(z)}{\zeta-z}\,dx\,dy.$$
Formally, $\bar\partial F(\zeta)=b(\zeta)$, so it suffices to check the convergence of the sum.  Noting that $|h_m(\zeta)|\ge\frac12$ in $U_m$ by Schwarz's lemma, termwise estimates give
    \begin{eqnarray*}
    \left|\frac1\pi\iint_{U_m}\frac{h_m(\zeta)}{h_m(z)}\frac{b(z)}{\zeta-z}\,dx\,dy\right|
        &\le& \frac2\pi|h_m(\zeta)|\iint_{U_m}\frac{|b(z)|}{|\zeta-z|}\,dx\,dy  \\
        &\le& \frac1\pi\Vert B\Vert_\infty|h_m(\zeta)|\frac{(1-|z_m|)^{-1}}{1-(4A^2)^{-1}}
            \iint_{B\left(z_m,\frac{1-|z_m|}{4A^2}\right)}\frac{dx\,dy}{|\zeta-z|}     \\
    \left|\frac1\pi\iint_{U_m}\frac{h_m(\zeta)}{h_m(z)}\frac{b(z)}{\zeta-z}\,dx\,dy\right|
        &\le& \frac4{3A^2}|h_m(\zeta)|\Vert B\Vert_\infty.
    \end{eqnarray*}
Summing, we find that
    $$\Vert F\Vert_\infty \le \frac{16}3\frac{\log4A^2}{\log\beta^{-1}}\Vert B\Vert_\infty$$
and
    $$|F(\zeta)|\le\frac{16}3\frac{\log4A^2}{\log\beta^{-1}}n^{-6}\Vert B\Vert_\infty
        \left(\frac{d(\zeta,\{\zeta_1,\zeta_2\})}{\diam(L_j)}\right)^2,$$
for $\zeta\in L_j$ due to the special properties of the functions $h_m$.\qed\\

Given this, the functions
    $$G_\mu(z)=\tilde g_\mu(z)+\sum_{\nu=1}^M(a_{\mu\nu}(z)-a_{\nu\mu}(z))f_\nu(z),$$
are corona solutions in $\mathbb D$ satisfying
    $$|G_\mu(z)-G_\mu^j(z)| \le \sum_{\nu=1}^M|f_\nu(z)||a_{\mu\nu}(z)-a_{\nu\mu}(z)|\lesssim M n^{-6}
        \left(\frac{d(z,\{\zeta_1,\zeta_2\})}{\diam(L_j)}\right)^2$$
on $L_j$.  For fixed $\tilde\alpha$, if the angle $\alpha$ defining the lozenges $L^n_{J,j}$ in (\ref{lozenges}) is sufficiently small then $\phi^n_J(L^n_{J,j})\subset L_j$ (reindexing as appropriate).  Mapping back to the domains $\tilde T^n_J$, since $(\phi^n_J)^{-1}$ behaves as $(z-\zeta_i(j))^{\frac{1+\alpha}2}$ about $\zeta_i(j)$, and so we obtain
    \begin{eqnarray*}
    |g^{(n,J)}_\mu(z)-g^j_\mu(z)| &\lesssim& n^{-6}\left(\frac{d(z,\{x_1,x_2\})}{\diam(L^n_{J,j})}\right)^{1+\alpha} \\
        &\lesssim& n^{-6}\frac{d(z,K)}{\diam(L^n_{J,j})},
    \end{eqnarray*}
on $L^n_{J,j}$, where $\{x_1,x_2\}=K\cap\overline{L^n_{J,j}}$, and $\Vert g^{(n,J)}_\mu\Vert_\infty\le C(M,\eta)$, with constants uniform in $(n,J)$.

This completes the proof of the proposition.\qed\\

\end{document}